\documentclass[a4paper,11pt]{article}
\usepackage{geometry,amssymb,amsmath,enumerate,latexsym,array}
\usepackage[dvips]{graphicx}
\usepackage[all]{xy}
\begin{document}
\parskip=1.5ex

\newtheorem{example}{Example}[section] 
\newtheorem{Def}[example]{Definition} 
\newtheorem{rem}[example]{Remark} 
\newtheorem{prop}[example]{Proposition} 
\newtheorem{cor}[example]{Corollary} 
\newtheorem{thm}[example]{Theorem} 
\newtheorem{lem}[example]{Lemma} 
\newenvironment{proof}{\noindent {\bf Proof} }{ \hfill 
$\Box$ \mbox{}\\} 
\newtheorem{alg}[example]{Algorithm}

 \title { Using Automata to obtain Regular Expressions for Induced Actions
 \thanks{ 
 \newline AMS 1991 CLASSIFICATION: 08A50, 16B50, 68Q40, 68Q42, 68Q45.  }}
\author{ Anne Heyworth \thanks{ Supported 1995-8 by an EPSRC
Earmarked Research Studentship, `Identities among relations for
monoids and categories', and 1998-9 by  a University of Wales, Bangor,
Research Assistantship.}\\ School of Mathematics \\ University of Wales,
Bangor \\
Gwynedd, LL57 1UT \\ United Kingdom \\ map130@bangor.ac.uk} 
\maketitle

\begin{center} Abstract \end{center}
{\small 
Presentations of Kan extensions of category actions provide a natural framework
for expressing induced actions, and therefore 
a range of different combinatorial problems. Rewrite systems for
Kan extensions have been defined and a variation on the Knuth-Bendix 
completion procedure can be used to complete them -- when possible.
Regular languages and automata are a useful way of expressing sets and actions, 
and in this paper we explain how to use rewrite systems for Kan extensions to
construct automata expressing the induced action and how sets of normal forms 
can be calculated by obtaining language equations from the automata. }

\section{Introduction}

Given a morphism of monoids $F: A \to B$ and an action of $A$ on a set $X$, 
the induced action of $B$ is on a set $F_*(X)$. Suppose $B$ has a 
presentation $mon\langle \Delta | RelB \rangle$ and $\Gamma$ is a set 
of generators for $A$ so that $F(a)$ is described in terms of $\Delta^*$,
and the action of $a$ on $X$ is known for each $a \in \Gamma$.
The problem is to describe $F_*(X)$. The usual rewrite theory is the case
where $A$ is the trivial monoid and $X$ is a one element set,
the extension to actions allows a wider range of applications. In fact our
extension goes beyond monoids to categories.

When $F: \mathsf{A} \to \mathsf{B}$ is a morphism of categories this gives a formulation 
in terms of induced actions of categories or
{\em Kan extensions}, as explained in \cite{paper2}, which 
defines rewrite systems for Kan extensions and introduces procedures for
completing such systems -- when possible.

This paper is a sequel to \cite{paper2}, showing how to interpret complete 
rewrite systems for Kan extensions. 
In this paper we assume that the completion procedure has been successful and
show how to use the rewrite systems to construct accepting automata whose 
languages can be calculated by equations, giving regular expressions for the
sets of the induced action (Theorem 4.3).
In the monoid case the induced action is on a single set. In the category
situation we may have many sets to describe. The use of languages is 
particularly appropriate for the situation where the action 
induced involves infinite sets.

Mac\,Lane  wrote that ``the notion of Kan extensions subsumes all the other 
fundamental concepts of category theory'' in section 10.7 of 
\cite{Mac} (entitled ``All Concepts are Kan Extensions'').
Together with \cite{paper2} this paper brings the power of rewriting 
theory and language theory to bear on a 
much wider range of combinatorial enumeration problems. 
Traditionally regular languages are used to specify the elements of a monoid
and rewriting is used for solving the word problem for monoids. 
Rewriting and regular languages may now also be used in the specification of 
\begin{enumerate}[i)] 
\item equivalence classes and equivariant equivalence classes, 
\item arrows of a category or groupoid, 
\item right congruence classes given by a relation on a monoid, 
\item orbits of an  action of a group or monoid. 
\item conjugacy classes of a group, 
\item coequalisers, pushouts and colimits of sets, 
\item induced permutation representations of a group or monoid. 
\end{enumerate} 
and many others.

\section{Rewrite Systems for Induced Actions}

This section gives a brief account of work of Brown and Heyworth
\cite{paper2} on extensions of rewriting methods.

Let $\mathsf{A}$ be a category. 
  A \textbf{category action} $X$ of $\mathsf{A}$ is a 
functor $X:\mathsf{A} \to \mathsf{Sets}$. 
Let $\mathsf{B}$ be a second category and let $F:\mathsf{A} \to \mathsf{B}$ be a  functor. 
Then an \textbf{extension of the action $X$ along $F$} is a pair $(K,\varepsilon)$ 
where $K:\mathsf{B} \to \mathsf{Sets}$ is a  functor and $\varepsilon:X \to K \circ F$ is 
a natural transformation.
The \textbf{Kan extension of the action $X$ along $F$} is an 
extension of the action $(K,\varepsilon)$ with the universal property that 
for any other extension of the action $(K',\varepsilon')$ there exists a 
unique natural transformation $\alpha:K \to K'$ such that 
$\varepsilon'=\alpha \circ \varepsilon$. 

The problem that has been introduced is that of ``computing a Kan 
extension''. Keeping the analogy with computation and 
rewriting for presentations of monoids and respecting the work of
\cite{BLW, CaWa1, CaWa2, Rosebrugh}, a definition of 
a presentation of a Kan extension is given as follows.

Recall that a \textbf{category presentation} is a pair $cat\langle \Delta | RelB \rangle$,
where $\Delta$ is a (directed) graph and $RelB$ is a set of relations on the free 
category $\mathsf{P}$ on $\Delta$. The category $\mathsf{B}$ presented by $cat\langle \Delta | RelB \rangle$
has objects $\mathrm{ Ob}\mathsf{B}$ that can be identified with $\mathrm{ Ob}\Delta$ and arrows $\mathrm{Arr}\mathsf{B}$
that can be identified with the classes of arrows of $\mathsf{P}$ under the congruence 
generated by $RelB$. The source and target functions of generating graph and category 
are denoted $src,tgt: \mathrm{Arr}\Delta \to \mathrm{ Ob}\Delta$ and $src,tgt:\mathrm{Arr}\mathsf{B} \to \mathrm{ Ob}\mathsf{B}$
respectively.

A \textbf{Kan extension data} $(X',F')$ consists of small categories
$\mathsf{A}$, $\mathsf{B}$ and functors $X':\mathsf{A} \to \mathsf{Sets}$ and $F':\mathsf{A} \to \mathsf{B}$.
A \textbf{Kan extension presentation} is a quintuple
$\mathcal{P}:=kan \langle \Gamma | \Delta | RelB | X | F \rangle$ where
\begin{itemize} 
  \item $\Gamma$ and $\Delta$ are (directed) graphs; 
  \item $X: \Gamma \to \mathsf{Sets}$ and $F: \Gamma \to \mathsf{P}$ are graph 
morphisms to the category of sets and the free category $\mathsf{P}$ on 
$\Delta$ respectively; 
  \item  and $RelB$ is a set of relations on the free category $\mathsf{P}$, i.e.
  a subset of $\mathrm{Arr}\mathsf{P} \times \mathrm{Arr}\mathsf{P}$. 
\end{itemize}

We say $\mathcal{P}$ \textbf{presents} the Kan extension $(K,\varepsilon)$ 
of the Kan extension
data $(X',F')$ where $X':\mathsf{A} \to \mathsf{Sets}$ and $F':\mathsf{A} \to \mathsf{B}$ if
\begin{itemize}
\item
 $\Gamma$ is a generating 
graph for $\mathsf{A}$  and $X: \Gamma \to \mathsf{Sets}$ is the restriction of
 $X':\mathsf{A} \to \mathsf{Sets}$ 
\item
$cat \langle \Delta | RelB \rangle$ is a category presentation for $\mathsf{B}$.
\item
$F:\Gamma \to \mathsf{P}$ induces $F':\mathsf{A} \to \mathsf{B}$.
\end{itemize}

We expect that a Kan extension $(K, \varepsilon)$ is given by 
\begin{itemize}
\item a set $KB$ for each $B \in \mathrm{ Ob} \Delta$, 
\item a function $Kb: KB_1 \to KB_2$ for 
each $b:B_1 \to B_2 \in \mathsf{B}$,  
\item a function $\varepsilon_A: XA \to KFA$ for each $A \in \mathrm{ Ob} \mathsf{A}$.
\end{itemize}

Let $\sqcup XA$  denote the disjoint union of all the sets $XA$ for all
objects $A$ in $\mathrm{ Ob}\mathsf{A}$ and let $\sqcup K\mathsf{B}$  denote the disjoint union of the
sets $KB$ for all objects $B$ in $\mathrm{ Ob}\Delta$.

 The main result of the paper \cite{paper2} defines rewriting procedures on the 
$\mathsf{P}$-set 
$$
T:= \bigsqcup_{B \in \mathrm{ Ob}\Delta} \bigsqcup_{A \in 
\mathrm{ Ob} \Gamma} XA \times \mathsf{P} (FA,B). 
$$ 
Elements of $T$ are called \emph{terms} and are written  $x|p$ where 
$x$ is an element of a set $XA$ for some object $A$ of $\mathsf{A}$, 
$p:FA \to B$ is an arrow of $\mathsf{P}$, and ``$|$'' is a symbol we use to separate 
the `element part' $x$ of the term from the `word part' $p$.

Then the set $T$ can also be written 
$$
T=\{x|p: x \in XA, p:FA \to B \text{ for some } A \in \mathrm{ Ob}\Gamma, B \in \mathrm{ Ob}\Delta\}.
$$

If $\mathcal{R}$ is a rewrite system on $T$ then we will write 
$\mathcal{R}=(\mathcal{R}_T,\mathcal{R}_P)$, since two kinds of rewriting are 
involved here. Rewriting using the rules $\mathcal{R}_P$ is the familiar 
$x|ulv \to x|urv$ given by a relation $(l,r)$. The rules
$\mathcal{R}_T$ derive from a given action of certain words 
on elements, so allowing rewriting $x|F(a) v \to x \cdot a|v$. 
Further, the elements $x$ and $x \cdot a$ may belong to different 
sets. When such rewriting procedures complete, the associated normal form 
gives in effect a computation of what we call the 
{\em Kan extension defined by the presentation}.

\begin{thm}[Data for Kan Extensions] \cite{paper2}\\ 
Let $\mathcal{P}=kan \langle \Gamma | \Delta | RelB | X | F \rangle$ be a 
Kan extension presentation. 

Let $\mathsf{P}$ be the free category on $\Delta$, let
$T:=\{x|p: x \in XA, p:FA \to B \text{ for some } A \in \mathrm{ Ob}\Gamma, B \in \mathrm{ Ob}\Delta\}$ 
and define $\mathcal{R}=(\mathcal{R}_\varepsilon, \mathcal{R}_K)$ where
$\mathcal{R}_\varepsilon:=
      \{(x|Fa,x\cdot a | i\!d_{FA_1}): x \in XA_1, a:A_1 \to A_2 \text{ in } \mathsf{A}\}$
and $\mathcal{R}_K:=RelB$. 

Then the Kan extension $(K,\varepsilon)$ presented by 
$\mathcal{P}$ may be  given by the following data: 
\begin{enumerate}[1)] 
\item 
the set $\sqcup KB = T/ \stackrel{*}{\leftrightarrow}_\mathcal{R}$, 
\item 
the function $\overline{\tau}: \sqcup KB \to \mathrm{ Ob} \mathsf{B}$ induced by 
$\tau:T \to \mathrm{ Ob} \mathsf{P}$, 
\item 
the action of $\mathsf{B}$ on $ \sqcup KB$ induced by the action of 
$\mathsf{P}$ on $T$, 
\item 
the natural transformation $\varepsilon$ determined by $x \mapsto [x| 
i\!d_{FA}]$ for $x \in XA, \;A \in \mathrm{ Ob} \mathsf{A}$. 
\end{enumerate} 
\end{thm}


To work with a rewrite system $\mathcal{R}$ on $T$ certain concepts of order on 
$T$ are required. The paper \cite{paper2} gives properties of  
orderings $>_X$ on $\sqcup XA$ and  $>_P$ on $\mathrm{Arr}\mathsf{P}$ which enable 
the construction of an ordering $>_T$ on $T$ with the properties needed for 
the rewriting procedures. For this paper we will assume that the order on
$T$ is a short-lexicographic order induced by ordering all the variables in
the alphabet we will be using.

Given a rewrite system $\mathcal{R}$ for a Kan extension and an ordering 
$>_T$ on $T$, a reduction relation $\to_\mathcal{R}$ compatible with the 
ordering is determined.
A reduction relation on a set is complete if it is Noetherian and confluent.
The Noetherian property implies that any term of $T$ can be repeatedly reduced 
until, after a finite number of reductions, an irreducible element will 
be obtained. The confluence property implies that if two terms are equivalent
under the relation $\stackrel{*}{\leftrightarrow}_\mathcal{R}$ then they reduce 
to the same term, i.e. there is a unique irreducible term in each equivalence 
class.
By standard abuse of notation the rewrite system $\mathcal{R}$ will be called 
complete when  is complete.  

The paper \cite{paper2} defines a variation on the Knuth-Bendix procedure which 
can be applied to $\mathcal{R}$ to complete it -- when this is possible. The procedure
has been implemented in {\sf GAP3} (to be converted to {\sf GAP4}), 
using a short-lex ordering.
The details in this paper show how to use automata to interpret the output of
the procedure when the sets $KB$ on which the induced action is defined cannot
be enumerated (i.e. are infinite).

\newpage
\section{Regular Languages and Automata for Induced Actions}
 
For a detailed introduction to automata theory, refer to \cite{Cohen} or 
\cite{Hopcroft}. This section only outlines the essential ideas we use.

A (finite) \textbf{deterministic automaton} is a 5-tuple $\underline{A}=(S, 
\Sigma, s_0,  
\delta, Q)$ \ where $S$ is a finite set of \emph{states} (represented by circles), 
$s_0 \in S$ is the \emph{initial state} (marked with an arrow), 
$\Sigma$ is a finite \emph{alphabet},  
$\delta:S \times \Sigma \to S$ is the \emph{transition}, 
$Q\subseteq S$ is the set of \textbf{terminal states} (represented by double 
circles).
A deterministic automaton $\underline{A}$ is \textbf{complete} 
if $\delta$ is a function, and \textbf{incomplete} if it is only a 
partial function. If $\underline{A}$ is incomplete, then when 
$\delta(s,a)$ is undefined, the automaton is said to \textbf{crash}.

The \textbf{extended state transition} $\delta^*$ is the extension of $\delta$ 
to $\Sigma^*$. It is defined by 
$\delta^*(s,i\!d):=s$, $\delta^*(s,a):=\delta(s,a)$, $\delta^*(s,aw) := 
\delta^*(\delta(s,a),w)$ where $s\in S$, $a\in \Sigma$ and $w\in \Sigma^*$.
We are interested in the final state $\delta^*(s_0,w)$ of the machine after a 
string $w \in \Sigma^*$ has been completely read. If the machine crashes 
or ends up at a non-terminal state then the string is said to have been 
\textbf{rejected}. If it ends up at a terminal state then we say the string 
is \textbf{accepted}.

A \textbf{language} over a given alphabet $\Sigma$ is a subset $L \subseteq 
\Sigma^*$. 
The set $L(\underline{A})$ of all acceptable strings is the 
 \textbf{language accepted by the automaton} $\underline{A}$. A language $L$ is a 
\textbf{recognisable} if it is accepted by some automaton $\underline{A}$. 
Two automata are \textbf{equivalent} if their languages are equal.
 
The \textbf{complement} of a complete, deterministic automaton is found by 
making non-terminal states terminal and vice versa. If the language accepted by 
an automaton $(\underline{A})$ is $L$, then the language accepted by its 
complement $(\underline{A})^C$ is $\Sigma^*-L$.

\begin{lem}[Completion of Automata] \cite{Cohen}\\
Let $\underline{A}=(S, \Sigma, s_0, \delta, Q)$ be an incomplete deterministic 
automaton. Then there exists a complete deterministic automaton 
$\underline{A}^{CP}$ such that $L(\underline{A})=L(\underline{A}^{CP})$. 
\end{lem} 

Diagrammatically this means that automata may be completed by adding one further 
non-terminal (dump) state $d$ and adding in all the missing arrows so that they 
point to this state.

A \textbf{non-deterministic automaton} is a 5-tuple 
$\underline{A}=(S, \Sigma, S_0, \delta, Q)$ \  
where $S$ is a finite set of states, 
$S_0 \subseteq S$ is a set of initial states, 
$\Sigma$ is a finite alphabet,  $Q\subseteq S$ is the set of 
terminal states and 
$\delta:S \times \Sigma \to \mathbb{P}(S)$ is the transition mapping 
where $\mathbb{P}(S)$ is the power set. 
The \textbf{language accepted by a non-deterministic automaton} $\underline{A}$
is the set of words $L \subseteq \Sigma^*$ such that $\delta^*(s,l) \cap Q \not=
\emptyset$ for some $s\in S$ for all $l \in L$.
 
\begin{lem}[Determinising Automata] \cite{Cohen}\\ 
Let $\underline{A}=(S, \Sigma, S_0, \delta_1, Q)$ be a non-deterministic 
automaton. Then there exists a deterministic automaton $\underline{A}^D$ such 
that $L(\underline{A}^D)=L(\underline{A})$. 
\end{lem} 

In practice a non-deterministic automaton may be made deterministic by drawing 
a \emph{transition tree} and then converting the tree into an automaton;  
for details of this see  \cite{Cohen}.

A regular expression is a string of symbols representing a regular language. 
Let $\Sigma$ be a set (alphabet). The empty word will be denoted $i\!d$.
A \textbf{regular expression} over $\Sigma$ is a string of symbols formed by 
the rules 
\begin{enumerate}[i)] 
\item $a_1\cdots a_n$ is regular for $a_1,\ldots,a_n \in \Sigma$, 
\item $\emptyset$ is regular, 
\item $i\!d$ is regular, 
\item if $x$ and $y$ are regular then $xy$ is regular, 
\item if $x$ and $y$ are regular then $x+y$ is regular, 
\item if $x$ is regular then $x^*$ is regular. 
\end{enumerate} 

For example $(x+y)^*-z$ is the expression representing the regular language 
$(\{x\} \cup \{y\})^* \, / \, \{z\}$. 
For our purposes a \textbf{right linear language equation} over $\Sigma$ is an 
expression $X=AX+E$ where $A,X,E \subseteq \Sigma^*$.
 
\begin{thm}[Arden's Theorem] \cite{Cohen}\\ 
Let $A,X,E \subseteq \Sigma^*$ such that $X=AX+E$ where $A$ and $E$ are known and $X$ is unknown. 
Then 
\begin{enumerate}[i)] 
\item $A^*E$ is a solution, 
\item if $Y$ is any solution then $A^*E \subseteq Y$, 
\item if $i\!d \not\in A$ then $A^*E$ is the unique solution. 
\end{enumerate} 
\end{thm} 
 
\begin{thm}[Solving Language Equations] \cite{Cohen}\\ 
A system of right linear language equations:
\vspace{-0.2cm} 
\begin{alignat*}{3} 
X_0 & = \, A_{0,0}X_0  &&+ \cdots + A_{0,n-1}X_{n-1} && +\, E_0\\ 
X_1 & = \, A_{1,0}X_0  &&+ \cdots + A_{1,n-1}X_{n-1} && +\, E_1\\ 
\cdots & \qquad \cdots &&\mbox{\hspace{4mm}} \cdots \qquad \cdots     && \quad \cdots\\ 
X_{n-1} & = \, A_{n-1,0}X_0 &&+ \cdots + A_{n-1,n-1}X_{n-1} && +\,  E_{n-1} 
\end{alignat*} 
where $A_{i,j},E_i \in \mathbb(\Sigma^*)$ and $i\!d \not\in A_{i,j}$ for 
$i,j=0,\ldots,n-1$, has a unique solution. 
\end{thm} 

\newpage
\begin{thm}[Regular Expressions from Automata] \cite{Cohen}\\ 
Let $\underline{A}$ be a deterministic automaton. 
Then $L(\underline{A})$ is regular. 
\end{thm} 

\begin{proof}
Let $\underline{A}=(S,\Sigma,s_0,\delta,Q)$, where $S:=\{s_0,\ldots,s_{n-1}\}$.
For $i=1,\ldots,n-1$ define $X_i:=\{z \in \Sigma^* : \delta^*(s_i,z) \in Q\}$.
It is clear that $X_0=L(\underline{A})$.

Define $E_i:= \{i\!d\}$ if $s_i \in Q$ and $\emptyset$ otherwise.

Define $A_{i,j}$ to be the sum of all letters $x \in \Sigma$ such that 
$\delta^*(s_i,x) = s_j$.

Then form the following system of equations:
\begin{alignat*}{3} 
X_0 & = \, A_{0,0}X_0  &&+ \cdots + A_{0,n-1}X_{n-1} && +\, E_0\\ 
X_1 & = \, A_{1,0}X_0  &&+ \cdots + A_{1,n-1}X_{n-1} && +\, E_1\\ 
\cdots & \qquad \cdots &&\mbox{\hspace{4mm}} \cdots \qquad \cdots     && \quad \cdots\\ 
X_{n-1} & = \, A_{n-1,0}X_0 &&+ \cdots + A_{n-1,n-1}X_{n-1} && +\,  E_{n-1} 
\end{alignat*} 

There are $n$ right linear equations in $n$ unknowns satisfying the conditions 
of Theorem 3.4. Therefore they have a unique solution.
\end{proof}

Thus every non-deterministic automaton gives rise 
to a system of language equations from whose solutions a description of the 
language may be obtained. 
 
\begin{thm}[Kleene's Theorem]\cite{Cohen}\\
A language $L$ is regular if and only if it is recognisable. 
\end{thm} 
 
This section has outlined the basic automata and language theory used in the 
paper.
Our main result (Theorem 4.3) is the construction, from a complete rewrite system for a Kan 
extension, of automata which recognise the elements of the extension as a 
regular language.

\section{Constructing and Interpreting the Automata}

Throughout this section we continue with the notation of \cite{paper2}
as described in Section 2.
Recall that a presentation of a Kan extension $(K,\varepsilon)$ 
is a 
quintuple $\mathcal{P}:=kan \langle \Gamma | \Delta | RelB | X | F \rangle$  
where $\Gamma$ and $\Delta$ are graphs, $RelB$ is a set of relations on 
the free category $\mathsf{P}$ on $\Delta$, while 
$X:\Gamma \to \mathsf{Sets}$ and $F:\Gamma \to \mathsf{P}$ are graph morphisms.
Recall that elements of the set
$$
T:=\bigsqcup_{B \in \mathrm{ Ob}\Delta} \bigsqcup_{A \in \mathrm{ Ob}\Gamma} XA \times \mathsf{P}(FA,B)
$$
are written $t=x|b_1 \cdots b_n$ with $x \in XA$, 
and $b_1,\ldots,b_n \in \mathrm{Arr}\Delta$ are composable with $src(b_1)=FA$.
The (`target') function $\tau:T \to \mathrm{ Ob}\Delta$ is defined by 
$\tau(x|b_1 \cdots b_n):=tgt(b_n)$ and the action of $\mathsf{P}$ on $T$, written
$t\cdot p$ for $t \in T$, $p\in \mathrm{Arr}\mathsf{P}$, is defined when $\tau(t)=src(p)$.

In \cite{paper2} we defined an initial rewrite system 
$\mathcal{R}_{init}:=(\mathcal{R}_\varepsilon,\mathcal{R}_K)$ on $T$ (also see Theorem 2.1),
and gave a procedure for attempting to complete this system. 
We will be assuming that the procedure has terminated, returning a complete 
rewrite system $\mathcal{R}=(\mathcal{R}_T,\mathcal{R}_P)$ with respect to a 
short-lex ordering on an alphabet $\Sigma$. 
In this section automata will be used to find regular expressions for 
each of the sets $KB$ for $B \in \mathrm{ Ob} \Delta$. 
 
Recall that $\sqcup XA$ is the union of the images under $X$ of all the objects 
of $\Gamma$ and $\sqcup KB$ is the union of the images under $K$ of all the 
objects of $\Delta$.
In general the automaton for the irreducible terms which are accepted as members 
of $\sqcup KB$ is the complement 
of the machine which accepts any string containing undefined compositions of 
arrows of $\mathsf{B}$, any string not containing a single $x_i$ on the left-most end, 
and any string containing the left-hand side of a rule. 
This essentially uses a semigroup presentation of the Kan extension. 
 
\begin{lem}[Semigroup presentation of a Kan extension]\mbox{ }\\ 
Let $\mathcal{P}$ present the Kan extension $(K, \varepsilon)$. 
Then the set $\sqcup KB$ may be identified with the non-zero elements of the 
semigroup having the presentation with generating set 
$$\Sigma_0:= (\sqcup XA) \sqcup \mathrm{Arr}\Delta \sqcup 0$$ 
and relations 
\begin{center} 
\begin{tabular}{ll} 
$0u=u0=0$     & for all $u \in \Sigma_0$,\\ 
$ux=0$        & for all $u \in \Sigma_0, \ x \in \sqcup XA$,\\  
$xb=0$        & for all $x \in XA, \ A \in \mathrm{ Ob}\Gamma, 
                         b \in \mathrm{Arr}\Delta$ such that $src(b)\not=FA$,\\ 
$b_1b_2=0$    & for all $b_1, b_2 \in \mathrm{Arr}\Delta$ such that 
$src(b_2)\not=tgt(b_1)$\\ 
$x(Fa)=(x\cdot a)$ & for all $x \in XA, \ a \in \mathrm{Arr} \mathsf{A}$ such that $src(a)=A$,\\ 
$l=r$         & for all $(l,r) \in RelB$.\\ 
\end{tabular} 
\end{center} 
\end{lem} 
 
\begin{proof} 
The semigroup defined is the set of equivalence classes of  
$T$ with respect to the second two relations (i.e. the Kan extension rules 
$\mathcal{R}_\varepsilon$ and $\mathcal{R}_K$) 
with a zero adjoined and multiplication of any two classes of $T$ defined to be 
zero. 
\end{proof} 
 
\begin{lem}[$T$ is a Regular Language]\mbox{ }\\
Let $\mathcal{P}$ be a presentation of a Kan extension $(K,\varepsilon)$. 
Then $T$ is a regular language over 
the alphabet $\Sigma:= (\sqcup XA) \sqcup \mathrm{Arr}\Delta$. 
\end{lem} 
\begin{proof} 
Define an automaton $\underline{A}:=(S,\Sigma,s_0,\delta,Q)$ where 
$S:=\{s_0, d\} \cup \mathrm{ Ob}\Delta$, $Q:=\mathrm{ Ob}\Delta$ and 
$\delta$ is defined as follows: 
\begin{align*}
\text{initial state } \quad \delta(s_0,u):=& 
\left\{ \begin{array}{ll} 
    FA & \text{ for } u \in XA, A \in \mathrm{ Ob} \Gamma\\ 
    d  & \text{ otherwise.}\\ 
        \end{array} \right. \\
\text{for } B \in \mathrm{ Ob} \Delta \quad
\delta(B,u):=&
\left\{ \begin{array} {ll} 
    tgt(u) & \text{ for } u \in \mathrm{Arr}\Delta, src(u)=B\\ 
    d & \text{ otherwise.}\\ 
                          \end{array} \right.\\ 
\text{dump state }\quad  \delta(d,u):= & \quad d \quad \text{ for all } u \in \Sigma. 
\end{align*} 
It is clear from the definitions that the extended state transition $\delta^*$ 
is such that $\delta^*(s_0,t) \in\mathrm{ Ob}\Delta$ if and only if $t \in T$. Hence 
$L(\underline{A})=T$. Therefore $T$ is regular over $\Sigma$. 
\end{proof}

\begin{thm}[Main Theorem] \mbox{ }\\
Let $\mathcal{R}$ be a finite complete rewrite system for the Kan extension
$(K,\varepsilon)$ given by the presentation 
$\mathcal{P}=kan\langle \Gamma|\Delta|RelB|X|F\rangle$. 

Let $T$ and $\Sigma$ be defined as before (Lemma 4.2). 
Then for each object $B \in \mathrm{ Ob}\Delta$  there is a regular expression 
representing a regular language $K_B$ over $\Sigma$ such that

\begin{enumerate}[i)] 
\item $KB=\{ [t]_\sim : t \in K_B \}$, where $[t]_{\mathcal{R}}$ represents the 
equivalence class of $t$ in $T$ with respect to 
$\stackrel{*}{\leftrightarrow}_{\mathcal{R}}$.
\item for $b: B_1 \to B_2$ in $\mathrm{Arr} \Delta$ the term $\mathsf{irr}(t \cdot b)$ is
an element of $K_{B_2}$
for all elements $t \in K_{B_1}$. 
\end{enumerate}

\end{thm} 
\begin{proof}
Recall the (`target') functions $tgt: \mathrm{Arr}\mathsf{P} \to \mathrm{ Ob}\Delta$ and 
$\tau:T \to \mathrm{ Ob}\Delta$. We use the following definition to restrict sets to 
those elements whose `target' is $B$.

$$H_B:=\left\{ \begin{array}{ll} 
    \{x : x \in XA \text{ and } XA=B, x \in H \} 
                & \text{ when } H \subseteq \sqcup XA\\ 
    \{p : tgt(p)=B, p \in H\}
                & \text{ when } H \subseteq \mathsf{P}\\
    \{t : \tau(t)=B, t \in H\}
                & \text{ when } H \subseteq T\\
         \end{array} \right.\\ $$       
                 
Then define
$\mathsf{irr}(H)$ be the set of irreducible forms of the terms $H \subseteq T$ with
respect to $\to_{\mathcal{R}}$.

For each object $B \in \mathrm{ Ob} \Delta$ we define an incomplete non-deterministic 
automaton $\underline{A}_B$ with input alphabet $\Sigma$, and language 
$\Sigma^*-\mathsf{irr}(T_B)$. This automaton rejects only the irreducible 
elements of $T_B$, i.e. it accepts all terms that do not represent elements of $T$,
terms that do not have `target' $B$ and terms that are reducible by 
$\to_{\mathcal{R}}$.

We will use the following notation:
\begin{center}
\begin{tabular}{l}
$\mathsf{l}(\mathcal{R}):=\{l:(l,r) \in \mathcal{R}\}$,\\
$\mathsf{pl}(\mathcal{R}):=\{u:(uv,r) \in \mathcal{R}\}$ and\\
$\mathsf{ppl}(\mathcal{R}):=\{u:(uv,r \in \mathcal{R}, v \not= i\!d \}$.
\end{tabular}
\end{center}
These are the set of all left hand side of rules, the set of all prefixes of 
left hand sides of rules and the set of all proper prefixes of left hand sides 
of rules respectively.

Now define $\underline{A}_B:=(S,\Sigma,s_0,\delta,Q_B)$ where

\begin{center} 
\begin{tabular}{l}
$S:=\{s_0,d\} \cup \mathrm{ Ob}\Delta \cup (\sqcup XA) \cup \mathsf{ppl}(\mathcal{R})$ and\\ 
$Q:=\{s_0,d, B\} \cup  (\sqcup XA)_B \cup\mathsf{ppl}(\mathcal{R})_B$.
\end{tabular} 
\end{center}

Let $x,b \in \Sigma$ so that $x \in \sqcup XA$ and $b \in \mathrm{Arr} \Delta$.
Define the transition 
$\delta: S \times \Sigma \to \mathbb{P}(S)$ by:
\begin{align*}
\text{initial state} \quad 
\delta(s_0,x):=&\left\{ \begin{array}{ll}
                  \{ x \}                & \text{ if } x \not\in \mathsf{l}(\mathcal{R}_T),\\
                  \{ d \}              & \text{ if } x \in \mathsf{l}(\mathcal{R}_T),
                  \end{array} \right.\\  
\delta(s_0,b):=& \{ d \},\\ 
\text{for }x_i \in XA \quad
\delta(x_i,x):=& \{ d \},\\
\delta(x_1,b):=&\left\{ \begin{array}{ll} 
                  \{x_1|b, tgt(b)\} & \text{ if } x_i|b \in \mathsf{ppl}(\mathcal{R}_T),\\ 
                  \{tgt(b)\}            & \text{ if } \tau(x_1)=src(b), x_i|b \not\in \mathsf{pl}(\mathcal{R}_T),\\ 
                  \{d\}                 & \text{ if } x_i|b \in \mathsf{l}(\mathcal{R}_T) \text{ or if } \tau(x_i) \not = src(b),\\
                  \end{array} \right.\\
\text{for }B_i \in \mathrm{ Ob}\mathsf{B} \quad  
\delta(B_i,x):=& \{d\},\\
\delta(B_i,b):=& \left\{ \begin{array}{ll} 
                  \{b, tgt(b) \}    & \text{ if } src(b)=B_i, b \in \mathsf{ppl}(\mathcal{R}_P),\\ 
                  \{tgt(b)\}            & \text{ if } src(b)=B_i, b \not\in \mathsf{pl}(\mathcal{R}_P),\\ 
                  \{d\}                 & \text{ if } src(b)=B_i, b \in \mathsf{l}(\mathcal{R}_P) \text{ or if } src(b) \not= B_i,\\
                   \end{array} \right.\\
\text{for }u \in \mathsf{ppl}(\mathcal{R}_T) \quad  
\delta(u,x):=& \{d\},\\
\delta(u,b):=& \left\{ \begin{array}{ll}
                  \{u \cdot b, tgt(b)\} & \text{ if } u \cdot b \in \mathsf{ppl}(\mathcal{R}_T),\\
                  \{tgt(b) \}               & \text{ if } \tau(u)=src(b), u \cdot b \not\in \mathsf{pl}(\mathcal{R}_T),\\
                  \{ d \}                     & \text{ if } u \cdot b \in \mathsf{l}(\mathcal{R}_T) \text{ or if } \tau(u) \not = src(b),\\
                   \end{array} \right.\\
\text{for }p \in \mathsf{ppl}(\mathcal{R}_P)  \quad
\delta(p,x):=& \{ d \},\\
\delta(p,b):=&\left\{ \begin{array}{ll} 
                  \{pb, tgt(b)\} & \text{ if } pb \in \mathsf{ppl}(\mathcal{R}_P),\\
                  \{tgt(b) \}         & \text{ if } tgt(p)=src(b), pb \not\in \mathsf{pl}(\mathcal{R}_P),\\
                  \{d\}              & \text{ if } pb \in \mathsf{l}(\mathcal{R}_P) \text{ or if } tgt(p) \not = src(b),\\
                  \end{array} \right. \\
\text{dump state} \quad
\delta(d,x):=& \{d\},\\
\delta(d,b):=& \{d\}.
\end{align*}
 
The extended state transition function $\delta^*$ is such that 
the intersection of $\delta^*(s_0,t)$ with  $Q_B$ is non-empty if and only
if $t$ is an element of $\Sigma^*$ which is not an element of $T_B$ or is 
reducible.\\

Thus for each object $B \in \mathrm{ Ob}\Delta$, and automaton $\underline{A}_B$ can be
constructed, where $L(\underline{A}_B)=\Sigma^*-\mathsf{irr}(T_B)$.
The results quoted in Section 3 allow us to make $\underline{A}_B$ deterministic
(Lemma 3.2) and take its complement. 
The language $K_B$ recognised by the resulting
automaton $(\underline{A}_B)^{DC}$ is $\Sigma^*-(\Sigma^*-\mathsf{irr}(T_B))$, i.e.
$K_B:=\mathsf{irr}(T_B)$. Hence (by Theorem 3.6) $K_B$ is regular.
Since $\mathcal{R}$ is a complete rewrite system on $T$ there exists a unique 
irreducible term in each class of $T_B$ with respect to 
$\stackrel{*}{\leftrightarrow}_{\mathcal{R}}$. Therefore the set $\mathsf{irr}(T_B)$ is
bijective with $T_B/ \stackrel{*}{\leftrightarrow}_{\mathcal{R}}=KB$.\\

The automaton $(\underline{A}_B)^{DC}$ gives rise to a system of right linear 
language equations (Theorem 3.5) with a unique solution, which is a regular 
expression for the language $K_B$ accepted by the automaton.
The regular expression can be obtained by applying Arden's Theorem (Theorem 3.3)
to solve the language equations.
Given that each set $KB$ is bijective with a regular language $K_B$, the action 
is described as follows: let 
$t \in K_{B_1}$ and $b:B_1 \to B_2$ for $B_1, B_2 \in \mathrm{ Ob} \Delta$, then 
$\mathsf{irr}(t \cdot b) \in K_{B_2}$.
\end{proof}

Thus for each object $B \in \mathrm{ Ob}\Delta$, an automaton $\underline{A}_B$ is
constructed, and a regular expression for the set $KB$ is obtained from solving
the language equations of the determinised complement of $\underline{A}_B$.
The $\mathsf{P}$-action on the elements $t$ of $T$ is right multiplication followed
by reduction with respect to $\to_{\mathcal{R}}$. This describes the functor 
$K$ in terms of regular expressions over $\Sigma^*$. The natural transformation 
$\varepsilon$ is given by $\varepsilon_A(x):=\mathsf{irr}(x|i\!d_{FA})$ for all $A \in \mathrm{ Ob} \mathsf{A}$
and $x \in XA$.

Therefore we have shown how the induced action $(K,\varepsilon)$ 
may be described in terms of regular languages and the 
reduction relation $\to_{\mathcal{R}}$.

\section{Example}

We construct simple automata which accept the terms which represent 
elements of some set $KB$ for $B \in \mathrm{ Ob}\mathsf{B}$ for an example of a Kan 
extension. The generating graphs are 
$$ 
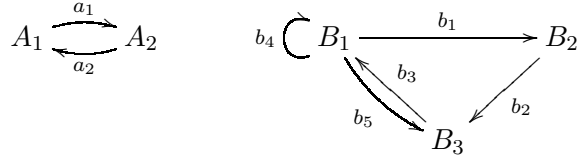
\xymatrix{A_1 \ar@/^/[r]^{a_1} 
        & A_2 \ar@/^/[l]^{a_2} 
       && B_1  \ar@(dl,ul)^{b_4} 
             \ar[rr]^{b_1} 
             \ar@/_/[dr]_{b_5} 
       && B_2 \ar[dl]^{b_2} \\ 
     &&&& B_3 \ar[ul]_{b_3} \\} 
$$ 
The relations are $RelB= \{ b_1b_2b_3=b_4 \}$, $X$ is defined by  
$XA_1 = \{ x_1, x_2, x_3 \}, XA_2 = \{ y_1, y_2 \}$   with 
$Xa_1:XA_1 \to XA_2: x_1 \mapsto y_1, x_2 \mapsto y_2, x_3 \mapsto 
y_1$, 
$Xa_2:XA_1 \to XA_2: y_1 \mapsto x_1, y_2 \mapsto x_2,$  \, and $F$ is 
defined by 
$FA_1=B_1$, $FA_2=B_2$, $Fa_1=b_1$ and $Fa_2 = b_2b_3$.\\ 
The completed rewrite system is: 
\begin{center} 
\begin{tabular}{llll} 
$x_1|b_1 \to y_1|i\!d_{B_2}$,      &  $x_2|b_1 \to y_2|i\!d_{B_2}$, 
  &  $x_3|b_1 \to y_1|i\!d_{B_2}$, &  $y_1|b_2b_3 \to x_1|i\!d_{B_1}$,\\ 
$y_2|b_2b_3 \to x_2|i\!d_{B_1}$,      &  $x_1|b_4 \to x_1|i\!d_{B_1}$, 
  &  $x_2|b_4 \to x_2|i\!d_{B_1}$,    &  $x_3|b_4 \to x_1|i\!d_{B_1}$,\\ 
$b_1b_2b_3 \to b_4$.             &&&\\ 
\end{tabular} 
\end{center} 
The proper prefix sets are 
$\mathsf{ppl}(\mathcal{R}_T):=\{y_1|b_2, y_2|b_2\}$ and 
$\mathsf{ppl}(\mathcal{R}_P):=\{b_1,b_1b_2\}$.
The following table defines the incomplete non-deterministic automaton 
which rejects only the terms of $T$ that are irreducible with respect to the 
completed relation $\to$. The alphabet over which the automaton is defined is 
$\Sigma:=\{x_1,x_2,x_3,y_1,y_2, b_1, b_2, b_3, b_4, b_5 \}$. 
\begin{center}
\begin{tabular}{l|llllllllll}
state/letter &$ x_1 $&$ x_2 $&$ x_3 $&$ y_1 $&$ y_2 $&$ b_1 $&$ b_2 $&$ b_3 $&$ b_4 $&$ b_5$\\
\hline
$s_0 $&$ x_1 $&$ x_2 $&$ x_3 $&$ y_1 $&$ y_2 $&$ d $&$ d $&$ d $&$ d $&$ d$\\
$x_1 $&$ d $&$ d $&$ d $&$ d $&$ d $&$ d $&$ d $&$ d $&$ d $&$ B_3$\\
$x_2 $&$ d $&$ d $&$ d $&$ d $&$ d $&$ d $&$ d $&$ d $&$ d $&$ B_3$\\
$x_3 $&$ d $&$ d $&$ d $&$ d $&$ d $&$ d $&$ d $&$ d $&$ d $&$ B_3$\\
$y_1 $&$ d $&$ d $&$ d $&$ d $&$ d $&$ d $&$ y_1|b_2, B_3 $&$ d $&$ d $&$ d$\\
$y_2 $&$ d $&$ d $&$ d $&$ d $&$ d $&$ d $&$ y_2|b_2, B_3 $&$ d $&$ d $&$ d$\\
$y_1|b_2 $&$ d $&$ d $&$ d $&$ d $&$ d $&$ d $&$ d $&$ d $&$ d $&$ d$\\
$y_2|b_2 $&$ d $&$ d $&$ d $&$ d $&$ d $&$ d $&$ d $&$ d $&$ d $&$ d$\\
$B_1 $&$ d $&$ d $&$ d $&$ d $&$ d $&$ b_1, B_2 $&$ d $&$ d $&$ B_1 $&$ B_3$\\
$B_2 $&$ d $&$ d $&$ d $&$ d $&$ d $&$ d $&$ B_3 $&$ d $&$ d $&$ d$\\
$B_3 $&$ d $&$ d $&$ d $&$ d $&$ d $&$ d $&$ d $&$ B_1 $&$ d $&$ d$\\
$b_1    $&$ d $&$ d $&$ d $&$ d $&$ d $&$ d $&$ b_1b_2, B_3 $&$ d $&$ d $&$ d$\\
$b_1b_2 $&$ d $&$ d $&$ d $&$ d $&$ d $&$ d $&$ d $&$ d $&$ d $&$ d$\\
$d      $&$ d $&$ d $&$ d $&$ d $&$ d $&$ d $&$ d $&$ d $&$ d $&$ d$\\
\end{tabular}
\end{center}

By constructing the transition tree for this automaton, 
we will make it deterministic. The next picture is of the partial transition 
tree -- the arrows to the node marked $\{d\}$ are omitted. 
$$ 
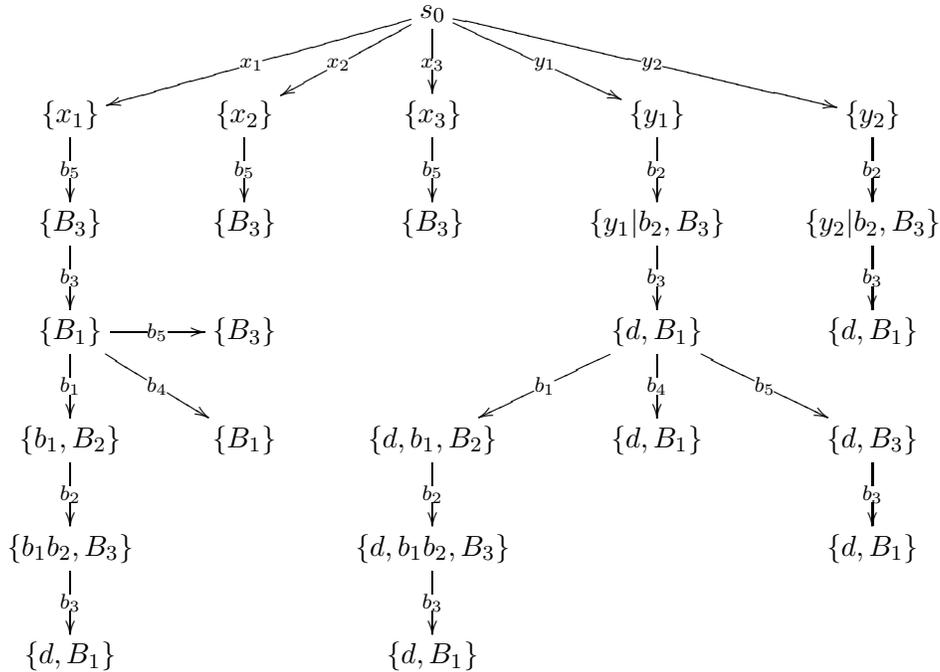
\xymatrix{&&s_0 \ar[dll]|{x_1} \ar[dl]|{x_2} \ar[d]|{x_3} \ar[dr]|{y_1} \ar[drr]|{y_2}&&\\
          \{x_1\} \ar[d]|{b_5} & \{x_2\} \ar[d]|{b_5} & \{x_3\} \ar[d]|{b_5} 
            & \{y_1\} \ar[d]|{b_2} & \{y_2\} \ar[d]|{b_2} \\
          \{B_3\} \ar[d]|{b_3} & \{B_3\} & \{B_3\}
            & \{y_1|b_2,B_3\} \ar[d]|{b_3} & \{y_2|b_2,B_3\} \ar[d]|{b_3}\\
          \{B_1\} \ar[d]|{b_1} \ar[dr]|{b_4} \ar[r]|{b_5} & \{B_3\} &
            & \{d,B_1\} \ar[dl]|{b_1} \ar[d]|{b_4} \ar[dr]|{b_5} & \{d,B_1\}\\
          \{b_1,B_2\} \ar[d]|{b_2} & \{B_1\} & \{d,b_1,B_2\} \ar[d]|{b_2}
            & \{d,B_1\} & \{d,B_3\} \ar[d]|{b_3}\\
          \{b_1b_2,B_3\} \ar[d]|{b_3} && \{d,b_1b_2,B_3\} \ar[d]|{b_3}
            & & \{d,B_1\}\\
          \{d,B_1\} && \{d,B_1\} &&\\}
$$ 
The tree is constructed with respect to the order 
on $\sqcup XA$ and $\mathrm{Arr}\Delta$, all arrows are drawn from $\{s_0\}$ and then
arrows from each new state created, in turn. 
When a label e.g. $\{B_3\}$ occurs that branch 
of the tree is continued only if that state has not been defined previously.
Eventually the stage is reached where no new states are defined, all the 
branches have ended. The tree is then converted into an automaton by `gluing' 
all states of the same label.
The initial state is  $\{s_0\}$ and a state is 
terminal if its label contains a terminal state from the original automaton. 
The automaton can often be made smaller, for example, here 
all the terminal states may be glued together.
One possibility is drawn below:\\
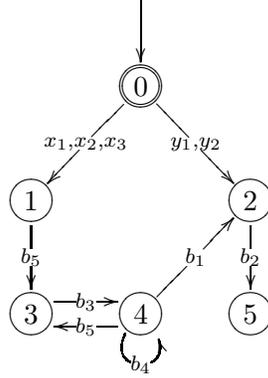
\vspace{-1cm} 
$$ 
\xymatrix{& \ar[d] & \\ 
         & *++[o][F=]{0} \ar[dl]|{x_1,x_2,x_3} \ar[dr]|{y_1,y_2} 
            & \\ 
            *++[o][F-]{1}  \ar[d]|{b_5}
         && *++[o][F-]{2}  \ar[d]|{b_2} \\
            *++[o][F-]{3}  \ar@<1ex>[r]|{b_3} &
            *++[o][F-]{4}  \ar@<1ex>[l]|{b_5} \ar@(dl,dr)|{b_4} \ar[ur]|{b_1}
            & *++[o][F-]{5} \\}
$$ 
Here the state labelled $1$, i.e. $S_1$ corresponds to the glueing together of 
$\{x_1\}$, $\{x_2\}$ and $\{x_3\}$ to form $\{x_1,x_2,x_3\}$ and the state $S_2$
is $\{y_1,y_2,b_1,B_2\}$.
States $S_3$ and $S_4$ represent $\{B_3\}$ and  $\{B_1\}$
respectively and state $S_5$ is $\{y_1|b_2,y_2|b_2,B_3,b_1b_2\}$.
The complement of this automaton accepts all irreducible elements of 
$\sqcup KB$. 
When $S_1$ and $S_4$ are terminal the language accepted is $K_{B_1}$. 
When $S_2$ is terminal the language accepted is $K_{B_2}$. 
When $S_3$ and $S_5$ are terminal the language accepted is $K_{B_3}$. 
The language equations from the automaton for $K_{B_1}$ are:
\vspace{-0.2cm} 
\begin{align*} 
X_0 & =(x_1+x_2+x_3)X_1+(y_1+y_2)X_2,\\ 
X_1 & =b_5X_3+i\!d_{B_1},\\ 
X_2 & =b_2X_5,\\ 
X_3 & =b_3X_4,\\ 
X_4 & =b_1X_2+b_4X_4+b_5X_3+i\!d_{B_1},\\ 
X_5 & =\emptyset.
\intertext{Putting $X_2=\emptyset$ and eliminating $X_1$ and $X_3$ 
by substitution gives}
X_0 &= (x_1+x_2+x_3)(b_5b_3X_4+i\!d_{B_1}),\\
X_4 &=(b_4+b_5b_3)X_4 + i\!d_{B_1}.
\intertext{Finally, applying Arden's Theorem to $X_4$
we obtain the regular expression}
X_0 &=(x_1+x_2+x_3)|(b_5b_3(b_4+b_5b_3)^*+i\!d_{B_1}).
\intertext{The separator ``$|$'' may be added at this point. 
Similarly, we can obtain regular expressions for $K_{B_2}$ and $K_{B_3}$. 
For $K_{B_2}$ we have}
X_0 &=(x_1+x_2+x_3)|b_5b_3(b_4+b_5b_3)^*b_1+(y_1+y_2)|i\!d_{B_2}.
\intertext{For $K_{B_3}$ the expression is}
X_0 &=(x_1+x_2+x_3)|(b_5b_3(b_4+b_5b_3)^*(b_1b_2+b_5) + b_5) + (y_1+y_2)|b_2.
\end{align*}

\end{document}